\documentclass[11pt,3p,times]{elsarticle}

\usepackage{enumitem}
\usepackage{hyperref}
\usepackage{amsmath}
\usepackage{pifont}
\usepackage[utf8]{inputenc}
\makeatletter
\def\LaTeX{\leavevmode L\raise.42ex
    \hbox{\kern-.3em\size{\sf@size}{0pt}\selectfont A}\kern-.15em\TeX}
\makeatother

\newcommand{\BibTeX}{{\rm B\kern-.05em{\sc
          i\kern-.025emb}\kern-.08em\TeX}}
\makeatletter
\def\@currentlabel{2.1}\label{e:dispaa}
\def\@currentlabel{2.21}\label{e:dispau}
\def\@currentlabel{2.22}\label{e:dispav}
\def\@currentlabel{2.23}\label{e:dispaw}
\def\@currentlabel{2.24}\label{e:dispax}
\def\theequation{\thesection.\@arabic\c@equation}
\makeatother

\renewcommand{\theequation}{\arabic{section}.\arabic{equation}}

\newcommand{\R}{\mathbb R}

\newcommand{\N}{\mathbb N}

\def \O{\Omega}
\everymath{\displaystyle}
\newtheorem{theorem}{Theorem}

\newtheorem{thm}{Theorem} [section]
\newtheorem{lem}{Lemma} [section]
\newtheorem{prop}{Proposition} [section]

\newtheorem{definition}{Definition} [section]

\newtheorem{rem}{Remark}[section]

\renewcommand{\theequation}{\thesection.\arabic{equation}}
\renewcommand{\thesection}{\arabic{section}}
\renewcommand{\theequation}{\thesection.\arabic{equation}}
\let\ssection=\section\renewcommand{\section}{\setcounter{equation}{0}\ssection}
\begin{document}
\begin{frontmatter}

\title{Existence and multiplicity of solutions for $m(x)-$polyharmonic elliptic Kirchhoff type equations without Ambrosetti-Rabinowitz conditions}
\author[mk0,mk1,mk2]{Mohamed Karim Hamdani\corref{cor1}}
\ead{hamdanikarim42@gmail.com}
\author[ah1,ah3,ah2]{Abdellaziz Harrabi}
\ead{abdellaziz.harrabi@yahoo.fr}
\cortext[cor1]{CONTACT\; Mohamed Karim Hamdani: hamdanikarim42@gmail.com}
\begin{center}
\address[mk0] {Science and technology for defense laboratory LR19DN01, center for military research, military academy, Tunisia.}
\address[mk1]{Military School of Aeronautical Specialities, Sfax, Tunisia.}
\address[mk2]{Mathematics Department, University of Sfax, Faculty of Science of Sfax, Sfax, Tunisia.}
\address[ah1]{ Mathematics Department, Northern Borders University, Arar, Saudi Arabia.}
\address[ah3]{Mathematics Department, University of Kairouan, Higher Institute of Applied Mathematics and Informatics, Kairouan, Tunisia.}
\address[ah2]{Senior associate in the Abdus Salam International Centre for Theoretical Physics, Trieste, Italy.}
\end{center}
\begin{abstract}
In this paper, we prove the existence of infinitely many solutions for a class of quasilinear
elliptic $m(x)$-polyharmonic Kirchhoff equations where the nonlinear function has a quasicritical growth at infinity and without assuming the Ambrosetti and Rabinowitz type condition. The new aspect consists in employing the notion of a Schauder basis to verify the geometry of the symmetric mountain pass theorem. Furthermore, for the case $m(x)\equiv Const$, we introduce a positive quantity $\lambda_M$ similar to the first eigenvalue of the $m$-polyharmonic operator to find a mountain pass solution, and also to discuss the sublinear case under large growth conditions at infinity and at zero. Our results are an improvement and generalization of the corresponding results obtained by Colasuonno-Pucci (Nonlinear Analysis: Theory, Methods and Applications, $2011$) and Bae-Kim (Mathematical Methods in the Applied Sciences, $2020$).
\end{abstract}

\end{frontmatter}
\section{Introduction}\label{nnn}
This paper is concerned with the existence and multiplicity of solutions to the following $m-$polyharmonic Kirchhoff problem:
\begin{eqnarray} \label{10c}
\begin{cases}
M\left(\int_{\O}|D_r u|^{m} \right)\Delta^r_m u = K(x)f(u) &\mbox{in}\quad \Omega, \\
u=\left(\frac{\partial}{\partial \nu}\right)^k u=0, \quad &\mbox{on}\quad \partial\Omega, \quad k=1, 2,.....  , r-1,
\end{cases}
\end{eqnarray}
where $\Omega \subset \mathbb{R}^N$ be a bounded smooth domain, $r \in \N^*$, $m >1$, $N\geq rm+1$, $K\in L^{\infty}(\O)$ is a positive weight function, $M \in C([0,+\infty))$ and $f\in
C(\mathbb{R})$ which will be specified later. The $m-$polyharmonic operator $\Delta^r_m$ is defined by

$$\Delta^r_m u=\begin{cases}-div \left\{\Delta^{j-1}(|\nabla \Delta^{j-1} u|^{m-2} \nabla \Delta^{j-1} u)\right\}, \text{ if $r=2j-1$}\\
\Delta^j (|\Delta^{j} u|^{m-2} \Delta^{j} u), \text{ if $r=2j$}
\end{cases}j\in \N^*,$$
which becomes the usual polyharmonic operator for $m=2$, namely $(-\Delta )^r$. Define the main $r-$order differential operator by
$$\textit{D}_r u=\begin{cases}
 \nabla\Delta^{j-1}u & \text{if $r=2j-1$},\\
\Delta^j u & \text{if $r=2j$},
\end{cases}\;\;j\in \N^*.
 $$
Note that $D_r u$ is an $N-$vectorial operator when $r$  is odd and $N > 1$, while it is a scalar operator when $r$ is even. Denote $E=W_0^{r,m}(\O)$   endowed with the following norm:
$$\|u\|^m=\int_{\O}|D_r u|^{m}.$$ $(E,\|\cdot\|)$ is a separable, uniformly convex, reflexive, real Banach space. Denote $E^*$ the  dual space of $E$, $p^* =\frac{mN}{N-rm}$ the Sobolev critical exponent and $q^*=\frac{p^*}{p^*-1}$ is the conjugate exponent of $p^*$. Recall also the Gagliardo-Nirenberg-Sobolev inequality which will be used a number of times later
\begin{eqnarray}\label{Gagliardo}
||u||_{L^{p^*}(\O) }\leq C\|u\|, \text{ for all } u \in E.
\end{eqnarray}

In recent years, there has been an increasing interest in studying problem \eqref{10c}, which has a broad background in many different applications, such as game theory, mathematical finance, continuum mechanics,
phase transition phenomena, population dynamics and minimal surface. The reader may consult \cite{BR,CP,CWL,HZ4,LLS1} and the references therein.

Problem \eqref{10c} is a general version of the following elliptic
problem with zero mass:
\begin{eqnarray*}\label{mm=1}
-\Delta u =K(x) f(u)\;\;
\mbox{ in }\O,
\end{eqnarray*}
for which some existence results have been established (see for example \cite{ASM,AZ}). In, Fan \cite{Fan}, studied a Kirchhoff type problem given by
 \begin{eqnarray}\label{lfaaa}
 -\Big(\varepsilon^2a+\varepsilon b\int_{\mathbb{R}^3}|\nabla u|^2dx\Big)\Delta u
  +V(x)u= K(x)|u|^{p-1}u,\quad  x\in \mathbb{R}^3,\nonumber\\
  u\in H^1(\mathbb{R}^3),
 \end{eqnarray}
 where $a,b$ are positive constants, $2< p<5$, $\varepsilon>0$ is a small parameter, $V(x)\in C^1(\mathbb{R}^3)$ and $K(x)\in C^1(\mathbb{R}^3,\mathbb{R})$ and
$0<K_\infty:=\lim_{|x|\to \infty}K(x)\leq K(x)$, and
$K(x)\not\equiv K_\infty$ for all $x\in\mathbb{R}^3$. With the help of Nehari-Pohozaev manifold, the authors proved the existence and concentration properties of a positive ground state solution as $\varepsilon\to 0$. Recently,  Li-Li-Shi \cite{LLS1} extended problem \eqref{lfaaa} to the following zero mass problem
\begin{eqnarray*}\label{prob-kirchh}
\left(a+\lambda\int_{\O}|\nabla u|^2\right)(-\Delta u)= K(x)f(u), \; \mbox{ in }\; \O,
 \end{eqnarray*}
where $\O=\R^N$, $\lambda\geq 0$  and under the assumptions
that $K$ is a nonnegative weight function satisfying
 $K \in [L^p(\O) \cap L^\infty(\O)]\setminus \{0\}$ for some $p\geq \frac{2N}{(N + 2)}.$ However, the authors used a cut-off functional and Pohozaev type identity to obtain the bounded Palais–Smale sequences.  Similar result has been obtained in \cite{LLS0} related to the delicate case $K(x)\equiv 1$, $N\geq3$, $a, b$ are positive constants and $\lambda\geq0$. 
 
The study of Kirchhoff type equations has already been extended to the case involving the $p$-Laplacian. In \cite{RFK}, Rasouli-Fani-Khademloo studied the following class of nonlocal Dirichlet problems
\begin{eqnarray}\label{clas}
\begin{cases}
-\left(a+b\int_{\R^N}|\nabla u|^pdx\right)\Delta_p u= K(x)f(u), \; \mbox{ in }\; \O,\\
u\in D^{1,p}(\O).
\end{cases}
 \end{eqnarray}
where $\O$ is a bounded domain in $\R^N$ , $N= 1, 2, 3$, $a, b$ are positive constants and $K(x)=1$. By using variational method and quantitative deformation lemma, they proved that problem \eqref{clas} possesses one sign-changing solution. Also, there are some recent works considered Kirchhoff type equations involving the $p$-Laplacian, see for instance \cite{Hamdani,HMH} and
the references therein.
\section{\bf Existence of multiple solutions for $\mathbf{\eqref{10c}}$}\label{SUBSECT1.1}
We assume that $K$ satisfies:\\
$(K_p)$:  $K \in L^\infty(\O)$ is a positive weight function.

We relax here the global structural assumption imposed on $M$ in \cite{CP} into the following:\\
$(M_1):$  $M: [0,+\infty) \rightarrow [0,+\infty)$ a continuous function satisfying: there exist $\tau_0\geq0$ and $\gamma\in (1,P_L=\frac{p^*}{m})$ such
 that $$\tau M(\tau)\leq \gamma \widehat{M}(\tau),\;\mbox{ for all } \tau\geq \tau_0, \mbox{ where } \widehat{M}(\tau)=\int_{0}^{\tau} M(z)dz.$$
Several recent papers for Kirchhoff problems deal only with the non-degenerate case, that is when $M(\tau)\geq s>0$ for all $\tau\in [0,+\infty) $ (see also \cite{CWL}). As in \cite{CP}, we covered the degenerate case at zero, we assume:\\\\
$(M_2):$  for each  $\eta>0 \mbox{ there exists } m_\eta>0 \mbox{ such that }
 M(\tau)\geq  m_\eta, \mbox{ for all } \tau\geq\eta.$

Set $F(s)=\int^s_0 f(t)$ and $p'=\frac{p}{p-1}$ the conjugate exponent of $p$. We introduce the following assumptions on the nonlinear term:
 \begin{enumerate}
\item[$(H_1):$]  there exists $C>0$ such that $C|f(s)|^{q^*}\leq s f(s)- m\gamma F(s)$,
 $\mbox{ for all } s\in\R;$
 \item[$(H_2):$]
$\lim\limits_{s \rightarrow \infty}$ $\frac{f(s)}{| s |^{p^*-1}} = 0$;
 \item [$(H_3)$:] $\lim\limits_{s\to
\infty} \dfrac{F( s)}{|s|^{m\gamma } }=\infty$.
\end{enumerate}

Problem \eqref{10c} is called nonlocal due to the
presence of the term
$ M\left(\int_{\O}|D_r u|^{m}\right),$
which implies that the equation in \eqref{10c} is no longer a pointwise identity. First, if we only assume that $M \in C([0,+\infty))$, $K \in L^\infty(\O)$ and $f \in C(\R)$ satisfying $(H_2)$, then the Euler-Lagrange functional associated to problem \eqref{10c} is defined by
\begin{eqnarray}\label{kir functi}
I(u)=\frac1m \widehat{M}(|D_r u|^{m})-\int_\O K(x)F(u),\;\mbox{ for all } u\in E
\end{eqnarray}
and $I \in C^1(E)$ with
\begin{eqnarray}\label{kir deriv}
\langle I'(u), {v}\rangle =M\left(\int_{\O}|D_r u|^{m}\right)\int_\O(|\textit{D}_r u|^{m-2}\textit{D}_r u\textit{D}_r v)  -\int_\Omega K(x)f(u)v,\; \mbox{ for all } u,v\in E.
\end{eqnarray}
So, $u \in E$ is a weak solution of \eqref{10c} if and only if $u $ is a critical point of $I$.

We begin by verifying the Palais-Smale condition (or the Cerami condition) under $(H_1)$-$(H_2)$ (see Remark \ref{r3} for further comments).

\begin{prop} \label{th13}
Assume that $(H_1)$-$(H_2)$ and $(M_1)$-$(M_2)$ hold. If, in addition $K$ satisfies $(K_p)$, then
\begin{itemize}
  \item[$1)$] $I$ satisfies the Palais-Smale condition
if $m\geq 2$;
  \item[$2)$] $I$ satisfies the Cerami condition if $1<m<
2$.
\end{itemize}
\end{prop}
We note that for $m \neq 2$, the variational setting of \eqref{10c} lacks an ordered Hilbert
space structure and so provokes some mathematical difficulties to obtain infinitely many solutions when $f(.)$ is an odd function. In a famous paper \cite{CP}, Colasuonno-Pucci established  multiplicity results by using minimax approach under a restrictive growth condition at zero (see also \cite{DJM}). Differently to \cite{CP,DJM} {\bf we will not here impose any control on $f$ at zero}. Also, we weaken the Ambrosetti-Rabinowitz condition assumed in \cite{BLW} by assumption $(H_1)$ (see Remark \ref{r3}) and we will not employ any cut-off technique as in \cite{BLW,LLS0,LLS1}. Our proof is a direct application of the following symmetric mountain pass theorem.
\begin{theorem} \label{th45}(\cite{R}).
 Let $E$  be a real infinite dimensional Banach space and $I\in C^1(E)$ satisfying the Palais-Smale condition; or the Cerami condition. Suppose $E=E^-\oplus E^+$, where
 $E^-$ is finite dimensional, and assume the following conditions:
\begin{enumerate}
  \item $I$ is even and $I(0)=0$;
 \item there exist $\alpha >0$ and $\rho>0$ such that $I(u)\geq \alpha$ for any $u\in E^+$ with $\|u\|=\rho$;
\item for any finite dimensional subspace $W\subset E$ there is $R= R(W)$ such that $I(u)\leq 0$ for $u\in W$, $\|u\|\geq R$;
\end{enumerate}
then, $I$ possesses an unbounded sequence of critical values.
\end{theorem}
The novelty here consists in using a Schauder basis of $E$ (see Corollary $3$ in \cite{FJN}) to verify point $2$ under only the condition $(H_2)$. Precisely, we have
\begin{lem} \label{lem113}
Assume that $(H_2)$ and $(M_2)$ hold. If, in addition $K$ satisfies $(K_p)$, then, {\bf for every} $\rho>0$ there exist a finite dimensional subspace $E^-$ and $\alpha >0$  such that $I(u)\geq \alpha$ for any $u\in E^+$ with $\|u\|=\rho$, where $E^+$ is a topological complement of $E^-$, i.e, $E=E^-\oplus E^+$.
\end{lem}
Our multiplicity existence result of this section reads as follows
\begin{thm}\label{th49}
Suppose that:
\begin{itemize}
 \item $M$ verifies $(M_1)$ and $(M_2)$;
  \item $f$ is an odd function verifying $(H_1)$-$(H_3)$;
  \item $K\in L^\infty(\O)$ is a positive weight function;
\end{itemize}
then, $I$ admits infinitely many distinct
pairs $(u_j, -u_j),\;j \in \N^*$, of critical points. Moreover,
$I(u_j)$ is unbounded.
\end{thm}
\section{\bf Existence results of general nonlinear elliptic equations}
In this section, we consider a general class of function that is $K(x)f(u)=f(x,u)$, and we rewrite \eqref{10c} in the following form:
\begin{eqnarray} \label{111}
\begin{cases}
M\left(\int_{\O}|D_r u|^{m} \right)\Delta^r_m u = f(x,u) &\mbox{in}\quad \Omega, \\
u=\left(\frac{\partial}{\partial \nu}\right)^k u=0, \quad &\mbox{on}\quad \partial\Omega, \quad k=1, 2,.....  , r-1.
\end{cases}
\end{eqnarray}
\subsection{\bf  Mountain pass solution:}
 To provide the mountain pass structure in the more familiar setting in the literature in which $M=1$, we require that $F(x,s)$ grows less rapidly than $ \frac{\lambda_1}{m}|s|^m$ near $0$ and more rapidly than $ \frac{\lambda_1}{m}|s|^m$ at infinity, where \begin{eqnarray}\label{lambda1}
\lambda_1:=\inf_{\substack{u\in E\\u\neq 0}}\dfrac{\int_\O|\textit{D}_ru|^m dx }{\int_\O|u|^m dx}>0,\end{eqnarray}
 is the first "eigenvalue" of $\Delta^r_m$. By analogy with $ \lambda_1$, we set
\begin{eqnarray}\label{lambdaM}\lambda_M :=\inf_{\substack{u\in E\\u\neq 0}}\dfrac{  \widehat{M}(\|u\|^m) }{\int_\O|u|^{m\gamma} dx},\end{eqnarray}
with $1<m\gamma <p^*$. To ensure that $\lambda_M $ is positive, we need the following  {\bf sufficient and necessary } coercivity condition: \\\\
$(M_3):$  there is a positive constant $C$ such that $C\tau^\gamma \leq \widehat{M}(\tau), \; \mbox{ for all } \tau \geq 0.$\\
\begin{lem}\label{eig.posi}
\begin{enumerate}
  \item[$(i)$] $\lambda_M  \mbox{ is positive {\bf if and only if} } M \mbox{ satisfies }(M_3);$\label{lam pos}
  \item[$(ii)$] if $M= C\tau^{\gamma-1}$, then $\lambda_M $ is attained.\label{resul2}
\end{enumerate}
\end{lem}
An instructive example will be given in the end of the proof of Lemma \ref{eig.posi}, where $\lambda_M$ is not attained.

To derive a mountain pass solution we assume that $(M_1)$ is global, i.e. $\tau_0=0$, we also replace $(H_3)$ by the following conditions at infinity and at zero:
\begin{enumerate}
\item[$(H'_3):$] $\limsup \limits_{s\to 0} \dfrac{F(x, s)}{|s|^{m\gamma }}< \dfrac{\lambda_{M}}{m } < \liminf\limits_{s\to
\infty} \dfrac{F(x, s)}{|s|^{m\gamma } }, \;{  uniformly \; in }\; \;\overline\Omega.$
 \end{enumerate}
Then, we have
\begin{thm} \label{th24}
Assume that $(H_1)$-$(H_2)$-$(H'_3)$, $(M_1)$ (with $\tau_0=0$) and $(M_3)$ hold.  Then, problem \eqref{111} has a nontrivial mountain pass solution.
\end{thm}
\begin{rem}
 Theorem \ref{th24} holds if we substitute $(M_3)$ by $(M_2)$ and $(H'_3)$ by the following strong condition:
$$\limsup \limits_{s\to 0} \dfrac{F(x, s)}{|s|^{m\gamma }}=0 \mbox{ and }\liminf\limits_{s\to
\infty} \dfrac{F(x, s)}{|s|^{m\gamma } }=\infty, \;{  uniformly \; in }\; \;\overline\Omega.$$
\end{rem}
\subsection {\bf The $m\gamma-$sublinear case}
 For $ 1<m\gamma <p^*$, we introduce the following $m\gamma-$sublinear growth condition at infinity:
\begin{enumerate}
   \item[$(H'_{1}):$] $\limsup \limits_{s\to \infty} \dfrac{f(x, s)}{|s|^{m\gamma-2 }s}< \lambda_{M}.$
\end{enumerate}
 Then, we have the following result.
\begin{prop} \label{sublinear}
 Assume that $f$ and $M$ verify respectively $(H'_1)$ and $(M_3)$, then
\begin{itemize}
  \item[$1)$] $I(u)\to \infty$ as $\|u\|\to \infty$.
  \item[$2)$] $I$ satisfies the Palais-Smale condition.
\end{itemize}
\end{prop}
If in addition we assume
\begin{enumerate}
\item[$(H'_2):$]
 $ \liminf\limits_{s\to
0} \dfrac{F(x, s)}{|s|^{m\gamma } }=\infty,\mbox{ uniformly in }\;\overline\Omega$,

 \end{enumerate}
and
\begin{enumerate}
\item [$(M'_2):$] \;$\widehat{M}(\tau)\leq \beta\tau^\gamma,$ for all $\tau\geq 0,$ where $\beta$ is a positive constant.
 \end{enumerate}
 Then, we applied the Krasnoselskii genus theory to establish the following multiplicity result
\begin{thm}\label{thm1.6}
  Assume that $f (x, .)$ is an odd function for all $x\in \O$ satisfying $ (H'_1)$ and $(H'_{2})$ and $M$ verifies $(M'_2)$ and $(M _3)$, then $I$ admits infinitely many distinct
pairs $(u_j, -u_j),\;j \in \N^*$, of critical points.
\end{thm}
Note that in Proposition \ref{sublinear} and Theorem \ref{thm1.6} we may substitute $ (H'_1)$ and $ (M_3)$ only by the following strong assumption
$$ \liminf\limits_{s\to
\infty} \dfrac{f(x, s)}{|s|^{m\gamma-2 }s }=0,\mbox{ uniformly in }\;\overline\Omega.$$
Next, if we relax $(H'_2)$ into the following condition:
\begin{enumerate}
\item[$(H''_2):$]$ \dfrac{\lambda_{M}}{m } < \liminf\limits_{s\to
0} \dfrac{F(x, s)}{|s|^{m\gamma } }, \;{  uniformly \; in }\;\overline\Omega$.
 \end{enumerate}
 Then, we have the following result.
\begin{thm}\label{th14}
 Assume that $ (H'_1)$, $(H''_{2})$ and $(M_3)$ hold. Then $I$ is bounded from below and $c=\inf\{ I(u),\; u \in E\}<0$ is a critical value of $I$. Consequently, problem \eqref{111} has a nontrivial solution.
\end{thm}
\section{A $m(x)-$polyharmonic Kirchhoff Equation}
In this section we extend the previous multiplicity results to a class of $m(x)-$polyharmonic Kirchhoff equations  in the following
form:
\begin{eqnarray}\label{high}
\begin{cases}
M\left( \int_\O \frac{|D_r u|^{m(x)}}{m(x)}dx\right)\Delta^r_{m(x)} u =f(x,u) \quad &\mbox{in}\quad \Omega, \\
 u=\left(\frac{\partial}{\partial \nu}\right)^k u=0, \quad &\mbox{on}\quad \partial\Omega, \quad k=1, 2,.....  , r-1.
\end{cases}
\end{eqnarray}
where, $\Delta^r_{m(x)}$ is called the $m(x)-$polyharmonic operator defined as
$$\Delta^r_{m(x)} u=\begin{cases}-\operatorname{div} \left\{\Delta^{j-1}(|\nabla \Delta^{j-1} u|^{m(x)-2} \nabla \Delta^{j-1} u)\right\}, \text{ if $r=2j-1$}\\
\Delta^j (|\Delta^{j} u|^{m(x)-2} \Delta^{j} u), \text{ if $r=2j$}
\end{cases}j\in \N^*.$$
As before, also here $\Omega$ is a bounded domain in $\R^N$ with the smooth boundary $\partial\O$. Throughout this paper, unless otherwise stated, we shall always assume that exponent $m\in C(\overline{\O})$ with
$$1<m_-:=\inf_\O m(x)\leq m(x)\leq m_+:=\sup_\O m(x)<\frac Nr$$
where $r$ is a given positive integer verifying $rm(x) < N$ and $p^*(x)$ denotes the critical variable exponent related to $m(x)$, defined for all $x\in\overline{\Omega}$ by the pointwise relation:
\begin{eqnarray*}
p^*(x) =\begin{cases}
           \frac{Nm(x)}{N-rm(x)}, & \mbox{if } rm(x)<N \\
           +\infty, & \mbox{if } rm(x)\geq N.
         \end{cases}
\end{eqnarray*}
In the following, we denote by $[W_0^{r,{m(x)}}(\O)]'$ the dual space of $W_0^{r,{m(x)}}(\O)$ and $q^*(x) = \frac{p^*(x)}{p^*(x)-1}$ the conjugate exponent of $p^*(x)$. The Kirchhoff function $M: [0,+\infty) \rightarrow [0,+\infty)$ is assumed to be continuous and to verify conditions $(M_1)$-$(M_2)$ given in Section \ref{SUBSECT1.1}, where here $P_L={(p^*)_-}/{m_+}$.
We assume that $f :\O\times\R\to \R$ is a continuous function, satisfying the following properties:
\begin{itemize}
\item[$(H^{exp}_1):$]  there exist $s_0>0$ and $C>0$ such that $C|f(x,s)|^{q^*(x)}\leq s f(x,s)- m_+\gamma F(x,s)$,
 $\mbox{ for all } |s|>s_{0}$ and $x\in \Omega;$
\item[$(H^{exp}_2):$] $\lim\limits_{s \rightarrow \infty}$ $\frac{f(x ,s)}{| s |^{p^*(x)-1}} = 0$ \quad uniformly in $x \in \Omega$;
\item[$(H^{exp}_{3}):$] $\lim\limits_{s\to
\infty} \dfrac{F(x, s)}{|s|^{m_+\gamma } }= \infty$ \quad uniformly in $ \overline\Omega;$
\item[$(H^{exp}_{4}):$] $f(x,-s)= -f(x,s)$ \quad for any $x \in \Omega, s \in\R$.
\end{itemize}

 We know that the extension from the $m$-polyharmonic operator to the $m(x)-$polyharmonic operator is not trivial, since the $m(x)$-polyharmonic operator possesses more complicated nonlinearities than the $m-$polyharmonic operator, mainly due to the fact that it is not homogeneous and some classical theories and methods, such as the Lagrange multiplier theorem and the theory of Sobolev space, cannot be applied. Recently, the variable exponent Lebesgue Sobolev spaces has been used in the last decades to model phenomena concerning nonhomogeneus materials, this is, a new field research and reflects a new type of physical phenomena, for example in nonlinear elasticity theory and in modelling electrorheological fluids (Acerbi and Mingione \cite{AM1}, Diening \cite{D}, Halsey \cite{H}, Ru\u{z}i\u{c}ka \cite{R3}, Rajagopal and Ru\u{z}i\u{c}ka \cite{RR2}) and from the study of electromagnetism and elastic mechanics (Hamdani and Repov\v{s} \cite{HR}, Zhikov \cite{Z}), and raise many difficult mathematical problems. After this pioneering models, many other applications of differential operators with variable exponents have appeared in a large range of fields, such as image restoration and image processing (Aboulaich et al. \cite{AMS08}, Li,Li and Pi \cite{LLP}) and mathematical biology (Fragnelli \cite{F}).

In the past years, there have been many contributions devoted to the higher order equation. The study of problems involving $m(x)-$biharmonic operators has been widely approached. For example, via  the mountain Pass Theorem, El Amrouss et al \cite{AMM} obtained the existence and multiplicity of solutions for  a class of $m(x)-$biharmonic equation of the form
\begin{equation}\label{biharmonic}
\begin{cases}
\Delta^2_{m(x)}u=\lambda |u|^{m(x)-2}u + f(x,u) \mbox{ in } \Omega, \\
u=\Delta u=0,\mbox{ on } \partial\Omega,
\end{cases}
\end{equation}
where $\O$ is a bounded domain in $\R^N$, with smooth boundary $\partial \O$, $\lambda\leq 0$ and $f$ satisfies the classical $(AR)$ condition. Afrouzi and Shokooh in \cite{Afrouz S}, studied the following Navier boundary-value problem
depending on two parameters and involving the $p(x)$-biharmonic
\begin{equation}\label{biharmonic af}
\begin{cases}
\Delta^2_{m(x)}u=\lambda f(x,u(x))+\mu g(x,u(x)) \mbox{ in } \Omega, \\
u=\Delta u=0,\mbox{ on } \partial\Omega,
\end{cases}
\end{equation}
where $\lambda$ is a positive parameter, $\mu$ is a non-negative parameter and  $f,g\in C^0(\Omega\times\mathbb{R})$. By using critical point theory, the authors proved the existence of infinitely many weak solutions for \eqref{biharmonic af}.

In \cite{HTT}, Heidarkhani et al., obtained in the case $\mu=0$ the existence of two solutions for the problem \eqref{biharmonic af} under some algebraic conditions with the classical $(AR)$ condition based on variational methods and critical point theory. Also for further studies on this subject, we refer the reader to \cite{AMC14,BR16,HHMR,KR2,LDP13,LT13,ZGC}.

 The existence and multiplicity of stationary higher order problems of Kirchhoff type (in $n$-dimensional domains, $n\geq 1$) have also treated in some recent papers, via variational methods like the Symmetric mountain pass theorem in \cite{new4,new3} and via a three critical point theorem in \cite{new2}. Especially, by covering the degenerate case, Colasuonno-Pucci established in \cite{CP} the existence of infinitely many solutions for problem \eqref{high} by using minimax approach. In  most of the aforementioned references with variable exponent, the standard ${(AR)}$ condition has played an important role which is originally due to Ambrosetti and Rabinowitz \cite{R}. Indeed, there are many superlinear functions which do not satisfy ${(AR)}$ condition. A simple computation proves that the following function
\begin{eqnarray}\label{yu}
f(x,s)=|s|^{m^+\gamma-2}s\ln^{q(x)}(|s|) \mbox{ with } q(x)> 1,
\end{eqnarray}
does not satisfy ${(AR)}$ but it is easy to see the above function \eqref{yu} satisfies $(H^{exp}_1)-(H^{exp}_3)$. For this reason, in recent years there were some authors studied the problem \eqref{high} trying to drop the condition ${(AR)}$. For instance, we refer the interested readers to  \cite{GWL,LRDZ} for the $m(x)$-Laplacian equation, \cite{ACM2,BR16,Zhou} for the $m(x)$-biharmonic equation. However, in literature the only results involving the $m(x)$-polyharmonic without assuming the ${(AR)}$ condition can be found in \cite{BK} which the authors extended the condition $(H_5)$  of \cite{LY} to the following assumption:
\begin{itemize}
  \item [$(F_4)$:] There exist $c_0 \geq 0$, $r_0 \geq 0$, and $k>\max\{1, \frac{N}{Lp_-}\}$ such that
  $$|F(x,s)|^k\leq c_0|s|^{kp_-}\mathfrak{F}(x,s)$$
\end{itemize}
for all $(x,s)\in\R^N\times \R$ and $|s|\geq r_0$, where $\mathfrak{F}(x,s)=\frac{1}{\theta p^+}F(x,s)s-F(x,s)\geq 0.$
\\
Clearly, our condition $(H^{exp}_1)$ is weaker also than $(F_4)$. For more details we refer the reader to point $5$ of Remark $1.1$ in the recent work of Harrabi-Hamdani-Selmi \cite{HHS} for the polyharmonic case.

Before stating our result, we give the definition of weak solutions for problem \eqref{high}:
\begin{definition}
We say that $u$ is a weak solution of problem \eqref{high} if $u$ satisfies
\begin{eqnarray*}
{M}\left( \int_\O \frac{|D_r u|^{m(x)}}{m(x)}dx\right)\int_\O|\textit{D}_r u|^{m(x)-2}\textit{D}_r u\textit{D}_r vdx=\int_\Omega f(x,u)vdx,
\end{eqnarray*}
for all $v\in W_0^{r,{m(x)}}(\O)$.
\end{definition}
In the light of the variational structure of \eqref{high}, we look for critical points of
the associated Euler-Lagrange functional $I:E\to \R$ defined as
\begin{eqnarray}\label{kir func}
I(u)=\widehat{M}\left( \int_\O \frac{|D_r u|^{m(x)}}{m(x)}dx\right)-\int_\O F(x,u)dx,\;\mbox{ for all } u\in E.
\end{eqnarray}
Note that $I$ is a $C^1(E, \R)$ functional and\begin{eqnarray}\label{kir deriv}
\langle I'(u), {v}\rangle ={M}\left( \int_\O \frac{|D_r u|^{m(x)}}{m(x)}dx\right)\int_\O|\textit{D}_r u|^{m(x)-2}\textit{D}_r u\textit{D}_r vdx  -\int_\Omega f(x,u)vdx,
\end{eqnarray}
for any $v\in W_0^{r,m(x)}(\O)$. Thus, critical points of $I$ are weak solutions of \eqref{high}.

The main result of this section is the following:
\begin{thm}\label{th46}
Assume that $f$ satisfies $(H^{exp}_{1})$-$(H^{exp}_{3})$ and $M$ verifies $(M_{1})-(M_{2})$, then, $I$ admits infinitely many distinct
pairs $(u_j, -u_j),\;j \in \N^*$, of critical points. Moreover,
$I(u_j)$ is unbounded.
\end{thm}
\section{Some Remarks}\label{r3}
In this section, we give some remarks and instructive examples in order to understand the improvement brought by our assumptions.
\begin{enumerate}
\item
  The two main assumptions that appeared in a rich literature ensuring the $(PS)$ condition  are the analogue of the Ambrosetti-Rabinowitz condition related to the Kirchhoff function $M$ (see \cite{CP, DJM, MMTZ, R}):
\begin{enumerate}
\item[${\bf(AR)_\gamma}:$] there are constants $\theta > m\gamma$ and $s_{0} > 0$ such that
 $ s f(x , s) \geq\theta F(x,s)>0,\; \mbox{ for all } |s| > s_{0} \mbox { and } x \in \Omega;$
\end{enumerate}
\begin{enumerate}
\item[${\bf(SCP)}:$]  there exist $C > 0$ and $p$ satisfying \;$\theta
 \leqslant p < p^*-1$ such that\; $| f(x , s)|  \leqslant C( | s |^{p } + 1), \; \mbox{ for all } (x ,s) \in \Omega  \times \mathbb{R}.$
\end{enumerate}
When $M\equiv1$ (and so $\gamma=1$), some attempts were made to relax conditions $(AR)_1$ and $(SCP)$ (see  \cite{BLW, CM, CSY, H, LL, LRRZ, LZ1, MS, SZ2} and the references therein). Point out that $(H_1)$-$(H_2)$ are weaker than $(AR)_\gamma$-$(SCP)$. In fact, $(SCP)$ implies $(H_2)$ and from  $(AR)_\gamma$ we have
 $0<(1-\frac {m\gamma}{\theta})sf(x,s)< sf(x,s) -m\gamma F(x,s),\; \mbox{ for all } |s|>s_0$ and $x\in \O.$ As $p^*-1=\frac{1}{q^*-1}$ then from $(H_2)$ we have   $|f(x,s)|^{q^*}\leq sf(x,s), \mbox{  for  } |s|\geq s'_0.$ Clearly, $(H_1)$ follows from the above inequalities.
\\
Note that, $(AR)_\gamma$ requires the following severe restriction called the strong $m\gamma$-superlinear condition:
\begin{enumerate}
  \item [${\bf(SSL)}$:] there exists $C>0$ such that\;\;
$F(x,s)\geq C |s|^\theta,\;\mbox{ for all } |s| > s_{0} \mbox { and }  x\in \Omega$.
\end{enumerate}
The most part of the literature used  conditions $(SCP)$ and $(SSL)$ to verify respectively points $2$ and $3$ of Theorem \ref{th45}  which are here weaken by $(H_2)$ and $(H_3)$.
\item Since we assume $\gamma<\frac{p^*}{m}$, then $(m\gamma-1)q^*\leq m\gamma$, so for $(m\gamma-1)q^*-1\leq\alpha<1$ and  $a > \gamma \lambda_M$, then
       a simple computation shows that $f_1(s)=a|s|^{m\gamma-2}s-|s|^{\alpha -1}s$ satisfies $(H_1)$ (and $(H_2)$-$(H'_3)$) but never $(AR)_\gamma$ nor $(SSL)$.
 \item In \cite{SZ2},  $(AR)_1$ was relaxed  into one of the following conditions (for $m=2$ and $M=1$): there are
constants $\theta > 2$ and $C > 0$ such that
 \begin{align}\label{bcz}
  | \theta F(x,s)-s f(x , s) |\leq C(1+s^2),\;\mbox{ for all } (x,s)\in \O\times \R,
  \end{align}
 or the {\bf global} convexity condition
 \begin{align}\label{bdz}
 H(x,s):= sf(x,s)-2F(x,s) \mbox{ is convex in } s,\;\mbox{ for all } x\in\Omega.
  \end{align} However, the following nonlinearity $f_2(s)=|s|^{m\gamma-2}s\ln^q(|s|)$ with $q\geq 1$ verifies
   $(H_1)$-$(H_3)$ and $(H'_3)$ but does not satisfies assumptions \eqref{bcz} and \eqref{bdz} if $q>1$.
\item  Let $ \gamma_1\in (1,\frac{p^*}{m}),\; \gamma_2 \geq \frac{p^*}{m}$. Consider the degenerate Kirchhoff function $M(\tau)=\tau^{\gamma_1-1} \mbox{ if } \tau\geq 1$ and
    $M(\tau)=\tau^{\gamma_2-1} \mbox{ if } \tau\leq 1$. We can see that $M$ satisfies $(M_1)$-$(M_2)$ but not the global assumption $(M)$ required in \cite{CP}.
\item
Consider the following Kirchhoff function introduced in \cite{CP}:
  $$ M(\tau)= a\tau^{\gamma_1-1} +b \tau^{\gamma_2-1} \mbox{ with } a\geq 0, b > 0 \mbox{ and } 1\leq \gamma_1\leq \gamma_2 <\frac{p^*}{m}.$$
Then, $M$ satisfies $(M_1)$-$(M_2)$ and also assumption $(M_3)$. Moreover $M$ is degenerate if $a=0$ or $a\neq0$ and $\gamma_1\neq 1$.\qed
 \end{enumerate}

The outline of this paper is the following: In Section \ref{section 6}, we give proofs of Proposition \ref{th13}, Lemma \ref{lem113} and Theorem \ref{th49}. Section \ref{section 7} is devoted to the proofs of Lemma \ref{eig.posi} and Theorem \ref{th24}. Also, we give the proofs of Proposition \ref{sublinear} and Theorems \ref{thm1.6}-\ref{th14} in the sublinear case. Finally, we treat in Section \ref{section 8} the multiplicity results for the more delicate case $m\equiv m(x)$ stated in Theorem \ref{th46}.

In the following, $|\cdot|$ denotes the Lebesgue measure in $\R^N$ and  $C$ (respectively $C_{\epsilon}$) denotes always a generic positive constant independent of $n$ and $\epsilon$ (respectively independent of $n$), even their value could be changed from one line to another one.  The integral $\int_\O u(x)dx$ is simply denoted by $\int_\O u(x)$.
\section{Proofs of Proposition \ref{th13}, Lemma \ref{lem113} and Theorem \ref{th49} \label{section 6}}
Let-us first establish some inequalities from assumptions $(M_1)$ and $(H_1)$-$(H_2)$ which will be useful to prove Proposition \ref{th13} and Lemma \ref{lem113}.

Set $p'=\frac {p}{p-1}>1$ the conjugate exponent of $p$.
 According to $(M_1)$ and $(H_1)$ , there exists $C_0>0$ such that
\begin{eqnarray}\label{er}
C|f(s)|^{q^*} \leq s f(s)- m\gamma F(s),\; \mbox{ for all } s\in  \R, \mbox{ and }
-C_0 \leq\gamma \widehat{M}(\tau)-\tau M(\tau), \; \mbox{ for all } \tau\geq 0.
\end{eqnarray}
  From $(H_2)$ it follows that for any $\epsilon>0$ there exists $C_{\epsilon}>0$ such that
\begin{eqnarray}\label{k10}
|f(s)|\leq \epsilon|s|^{p^*-1}+C_{\epsilon}
\end{eqnarray}
and
\begin{eqnarray}\label{pp}
|F(s)|\leq \epsilon|s|^{p^*}+C_{\epsilon}.
\end{eqnarray}
 Multiplying \eqref{k10} by $|s-s'|, s'\in \R$, and applying Young's inequality, then we derive
 \begin{eqnarray}\label{llvw}
|f(s) (s-s')|\leq\epsilon(|s|^{p^*}+ |s-s'|^{p^*}) +C_{\epsilon},\;\mbox{ for all } s\in \R.
\end{eqnarray}
We also recall some known results which will be essential to prove Proposition \ref{th13}. Consider the functional $\psi(u)= \frac{1}{m}\|u\|^{m},\; u\in E$, we have $ \psi \in C^1(E)$ with Fr\'echet's derivative
$$\langle \psi'(u), v\rangle=\int_{\R^N}|\textit{D}_r u|^{m-2}\textit{D}_r u\textit{D}_r v, \mbox{ for all } u,v \in E.$$
Set $\varphi(t)=t^{m-1}, t\geq 0$, clearly we have $\langle \psi'(u), u\rangle=\varphi(\|u\|)\|u\|$ and  it follows from H\"older's inequality that $\|\psi'(u)\|_{E^*}=\varphi(\|u\|)$.
Obviously, $\varphi$ is a normalization function and since $E$ is locally uniformly convex and so uniformly convex and reflexive Banach space, then the corresponding duality mapping $J_\varphi$ is single valued (i.e., $J_\varphi=\psi'$) and satisfies the $S_+$ condition (see Proposition $2$ in \cite{DJM}, respectively Lemma $3.2$ in \cite{LZ0}):
\begin{eqnarray}\label{khc}
\mbox{ if } u_n\rightharpoonup u \mbox{ and }\limsup_{n\to +\infty}\psi'(u_n)(u_n-u)\leq 0,\;\mbox{ then } u_n \rightarrow u.
\end{eqnarray}
\subsection{ \bf Proof of Proposition \ref{th13}.}
 Since we assume that $N> rm$, we may easily see that
 \begin{eqnarray}\label{err.}
 m > \frac{1}{q^*}+1 \mbox{ if  } m\geq 2.
\end{eqnarray}
Let $u_n$ be a $(PS)$ sequence of $I_{K,f}$ if $m\geq 2$ ({ respectively} $(C)$ sequence if $1< m< 2$). The case $u_n$ admits a subsequence which converges strongly to $0$ in $E$ is trivial. Hence, we may suppose that there exist $\eta_0>0$ and $n_0 \in \N$ such that $\|u_n\|^m\geq \eta_0$ for all $n\geq n_0$. So, in view of $(M_2)$ we can find $m_{\eta_0}>0$ such that
 \begin{eqnarray}\label{ing1.}
 m_{\eta_0}\|u_n\|^m \leq M(\|u_n\|^m)\|u_n\|^m, \forall n\geq n_0.
\end{eqnarray}
\\
\textbf{Step 1.} We shall prove that $u_n$ is bounded in $E$. First, from \eqref{kir deriv} and  \eqref{ing1}, we have
\begin{align*}
m_{\eta_0}\|u_n\|^m \leq  &\langle I'(u_n),u_n\rangle+\int_{\O} Kf(u_n)u_n.
\end{align*}
Apply H\"older's inequality to the second term in the right-hand side and using \eqref{Gagliardo}, we obtain
\begin{align}\label{kirch hk}
 m_{\eta_0}\|u_n\|^m \leq \langle I'(u_n),u_n\rangle&+ C\left(\int_{\O} K|f(u_n)|^{q^*}\right)^\frac{1}{q^*}\|u_n\|.\end{align}
From \eqref{kir func} and \eqref{kir deriv}, one has
\begin{eqnarray}\label{ing2}
\int_\O K \left [ f(u_n)u_n-m\gamma F(u_n)\right]= \left[\gamma
\widehat{M}(\|u_n\|^m)-M(\|u_n\|^m)\|u_n\|^m\right]+ \langle I'(u_n), {u_n}\rangle-m\gamma I(u_n).
 \end{eqnarray}
Taking into account that $u_n$ is a $(PS)$ sequence if $m\geq 2$ ( respectively $(C)$ sequence if $1< m< 2$), then from \eqref{er} we deduce
\begin{equation*}\int_{\Omega} K| f( u_{n}) |^{q^*}\leq C(1+\|u_{n}\| ), (\mbox{  respectively }\int_{\Omega}K | f( u_{n})
|^{q^*}\leq  C).
  \end{equation*}
 Combining now the above inequality with \eqref{kirch hk}, it follows
\begin{equation*}
\|u_{n}\|^{m}\leq C( 1 + || u_{n} ||^{{\frac{1}{q^*}+1}} )\mbox{ if } m\geq 2, (\mbox{  respectively } \|u_{n}\|^{m}
\leq C \| u_{n} \| \; \mbox{ if } 1<m<2).
\end{equation*}
 Therefore, thanks to \eqref{err.} the $(PS)$ sequence is bounded in $E$ if $m\geq 2$, and clearly, the $(C)$ sequence $u_n$ is also bounded if $1<m<2$.
 \\\\
{\bf Step 2.} We shall prove that the bounded sequence $u_n$ has
a strong convergent subsequence in $E$. In fact, we can find a subsequence (denoted again by $u_n$) and $u \in E$ such that $u_n$ converges to
$u$ weakly in $E$, \text{ and a.e. in } $\O$. Also $u_n$ and $u_n - u$ are bounded in $L^{p^*}(\O)$. Using again \eqref{kir deriv}, we get
\begin{eqnarray}\label{e:1221zz.}
M(\|u_n\|^m)\int_\O\left(|\textit{D}_r u_n|^{m-2}\textit{D}_r u_n\textit{D}_r (u_n-u)\right) =\langle I'(u_n), u_n-u\rangle+\int_\Omega K f(u_n)(u_n-u).
\end{eqnarray}
First we claim that
 \begin{eqnarray}\label{ggf}
                     \int_\Omega Kf(u_n)(u_n-u)\mbox{ converges to } 0.
                     \end{eqnarray}
So, we may assume that $u_n$ converges to $u$ in $L^1(\O)$.
Apply inequality \eqref{llvw} (with $s= u_n$ and $s'=u$, and integrate over $\O$, we obtain
 $$\left|\int_\Omega K f(u_n)(u_n-u)  \right|\leq \epsilon + C_\epsilon \int_{\O}|u_n-u|.$$
As $u_n$ converges strongly in $L^1(\O)$, then \eqref{ggf} follows.

Next, since $I'(u_n) \to 0$ in $E^*$ and $(u_n-u)$ is bounded in $E$, then from \eqref{e:1221zz.} and \eqref{ggf} we can deduce that $$M(\|u_n\|^m)\int_\O\left(|\textit{D}_r u_n|^{m-2}\textit{D}_r u_n\textit{D}_r (u_n-u)\right) \to 0.$$
Hence,  \eqref{ing1.} yields  $\int_\O\left(|\textit{D}_r u_n|^{m-2}\textit{D}_r u_n\textit{D}_r (u_n-u)\right)  \mbox{ converges to } 0.$ Invoking now the $S_+$ property (see \eqref{khc}), we conclude that $u_n$ converges strongly to $u$ in $E$. \qed
\subsection{ \bf  Proof of Lemma \ref{lem113}.}
 Let $(e_i)_{i \in \N^*}$ be a Schauder basis of $E$ (see Corollary $3$ in \cite{FJN} and also \cite{FHHSPZ,T}), which means that each $x\in E$ has a {\bf
unique representation} $x=\sum_{i=1}^{\infty}a_ie_i$, where $a_i$ are real numbers. Set $ E_j=span(e_1, e_2,..,e_j)$, then, the linear projection onto $ E_j$ i.e.,
$P_j: E\to E_j,\; P_j(x)= \sum_{i=1}^{j}a_ie_i$ is a {\bf continuous} linear operator for all $j\in \N^*$ (see \cite{S2}) \footnote{More precisely, $P_j$ are
 uniformly bounded, that is there exists $C>0$ such that $\|P_j(x)\|\leq C\|x\| $ for each $j\in \N^*$ and all $x\in E$ (see \cite{FHHSPZ,S2})}.
  Therefore, $F_j= N( P_j)$ (the kernel of $P_j$) is a topological complement of $E_j$, that is $E=E_j\oplus F_j$.\\
    Fix $\rho\geq 0$ and set $S_j^\perp(\rho)=\left\{u\in F_j \mbox{ such that }\|u\|=\rho\right\}$ and $\beta_j:=\sup_{S_j^\perp(\rho)} \int_{\Omega} K|F(u)|$.\\
  So, we claim that
 \begin{eqnarray}\label{uxy}
 \beta_j \rightarrow 0, \mbox{ as } j\rightarrow +\infty,
 \end{eqnarray}
  Before proving the claim \eqref{uxy}, let-us first end the proof of Lemma \ref{lem113}. Fix $\rho >0$, then for all $u\in S_j^\perp(\rho)=\left\{u\in F_j \mbox{ such that }\|u\|=\rho\right\}$, we have
$$I(u)\geq \frac 1m \widehat{M}( ||u||^m)-\int_\O K F(u)\geq \frac 1m \widehat{M}( \rho^m) -\beta_j.$$
Set $\alpha= \frac {1}{2m} \widehat{M}(  \rho^m)$, since we assume that $M$ satisfies $(M_2)$, then $\alpha >0$.
As $\beta_j$ converges to $0$, we can choose $j=j_0$ large enough such that $\beta_{j_0} \leq  \alpha $. Hence,
 $I(u)\geq  \alpha$, and so Lemma \ref{lem113} holds with $E^-=E_{j_0}$, $E^+=F_{j_0}$ and $\alpha= \frac {1}{2m} \widehat{M}(  \rho^m)$.
\\\\
{\bf Proof of \eqref{uxy}.}
 We argue by contradiction. Suppose that there exist $m_0>0$ and a subsequence (denoted by $\beta_j$) such that\\
$$m_0<\beta_j,\forall j\in \N^*.$$
From the definition of $\beta_j$, there exists $u_j\in S_j^\perp(\rho)$ such that
\begin{eqnarray}\label{uuk..}
m_0<\int_{\O} K|F(u_j)|\leq \beta_j.
\end{eqnarray}
   As $\|u_j\|=\rho$, then there exist a subsequence (denoted by $u_{j}$) and $u\in E $ such that $u_{j}$  converges weakly to $u$ and $\mbox{ a.e in }\; \Omega$, also $u_j$ is bounded in $L^{p^*}(\O)$. Fix $k\in \N^*$, as $P_k\circ P_{k+1}=P_k$, then $ F_{j}\subset F_k$ for all $j\geq k$, and so
 \begin{eqnarray}\label{uuxk}
u_{j} \in F_{k},\mbox{ that is } P_k(u_{j})\; \forall j\geq k.
\end{eqnarray}

 On the other hand, Since $P_k$ is a continuous linear operator and $u_{j}$  converges weakly to $u$, then $P_k(u_{j})$
 converges weakly to $P_k(u)$ which with \eqref{uuxk} implies that $P_k(u)=0$ and therefore $u=0$ as $u=\lim_{k\to
\infty}\sum_{i=1}^{k}a_ie_i=\lim_{k\to
\infty}P_k(u)$. Consequently, $u_j$ converges to $0$ $\mbox{ a.e in } \Omega$, and then $KF(u_{j})  \rightarrow 0 \mbox{ a.e in }\; \Omega$ as $KF(0)=0$. Substitute $s$ by $u_j$ in inequality \eqref{pp} and integrate over a measurable set $A\subset\O$, then if $\O$ is a bounded domain we derive
\begin{eqnarray*}
\int _{A} K|F(u_{j}) | \leq   C\epsilon  \int _{A} | u_{j}|^{p^*}
+ C_{\epsilon}|A| \leq C\epsilon.
\end{eqnarray*}
Hence if  $|A|<\frac{\epsilon}{C_{\epsilon}}$ (where $|A|$ denotes the Lebesgue measure of $A$), we deduce
\begin{eqnarray*}
\int _{A} K|F(u_{j}) |  \leq C'\epsilon.
\end{eqnarray*}
Taking into account that $\O$ is a bounded domain, then Vitali's theorem implies that $K|F(u_{j})|\rightarrow 0$ in $L^{1}(\O)$ and in view of \eqref{uuk..}, we obtain
$$
0< m_0\leq 0.
$$
 Thus, we reach a contradiction.\qed

\subsection{ \bf Proof of Theorems \ref{th49}}
We will show that the functional $I$ satisfies all conditions of the abstract Theorem \ref{th45}. In fact, Since $f$ is odd and $F(0)=0$, then $I$ is an
even functional and $I(0)=0$, and according to Proposition \ref{th13}, $I$ satisfies the $(PS)$ condition if $m\geq 2$ (respectively the $(C)$ condition if $1<m<2$). Thanks to Lemma \ref{lem113} $I$ verifies point $2$ of Theorem \ref{th45}. Therefore, it remains to show that condition $3$ of Theorem \ref{th45} holds. In view of $(H_3)$ we have for all $A>0$, there exists $C_A>0$ such that $$F(x,s)\geq A|s|^{m\gamma}- C_A,\; \forall (x,s)\in
\O\times \mathbb{R}.$$ Therefore we may conclude $$ I(u) \leq
\frac{C_1}{m}||u||^{m\gamma}-A\|u\|^{m\gamma}_{L^{m\gamma}(\O)} +C_A,\;\mbox{ for all }\|u\|>1.$$
The desired result followed which completes the proof of Theorem \ref{th49}.\qed
\section{Proofs of  Lemma \ref{eig.posi}, Proposition \ref{sublinear} and Theorems \ref{th24},\ref{thm1.6},\ref{th14}} \label{section 7}
\subsection{ \bf Proof of Lemma \ref{eig.posi}} { \bf Proof of $(i)$:} Assume that $M$ satisfies $(M_3)$. As $1<m\gamma <p^*$, Sobolev's inequality implies
$\int_\O|u|^{m\gamma} \leq C \|u\|^{m\gamma},$ which combined with $(M_3)$ yields  $\lambda_M >0$.

 Conversely, if $M$ does not satisfies $(M_3)$, then there is a sequence $\tau_i > 0$ such that $\tau_i^{-\gamma}\widehat M(\tau_i) \to 0$. Consider $\varphi \in E$ such that $\|\varphi\| = 1$, set $u_i = \tau_i^{1/m}\varphi$, then $\|u_i\|^m = \tau_i$. So
$$\frac{\widehat M(\|u_i\|^m)}{\int_\O|u_i|^{m\gamma}} = \frac{\widehat M(\tau_i)}{\tau_i^\gamma} \times \frac{1}{\int_\O|\varphi|^{m\gamma} } \to 0.$$ Therefore, $\lambda_M = 0$.\qed\\\\
{ \bf Proof of $(ii)$:} $(a)$ We will first prove that if $M= C\tau^{\gamma-1}$, then $\lambda_M$ is attained. without losing
any generality, we may assume that $C=\gamma$.
Let
\begin{eqnarray*}\label{lambdaM}\lambda_M :=\inf_{\substack{u\in E\\u\neq 0}}\dfrac{  \widehat{M}(\|u\|^m) }{\int_\O|u|^{m\gamma} }=\inf_{\substack{u\in E\\u\neq 0}}\dfrac{ \frac{C}{\gamma}\|u\|^{m\gamma}}{\|u|_{L^{m\gamma}(\O)}^{m\gamma} }=\inf\left\{\|u\|^{m\gamma}, \mbox{ such that } \|u|_{L^{m\gamma}(\O)}=1\right\}.\end{eqnarray*}

Let $u_n$ be a minimizing sequence, i.e, \; $\|u_n\|_{L^{m\gamma}(\O)}=1$ and \;$\|u_n\|^{m\gamma}\to \lambda_M$, so $u_n$ is bounded in the $E$ norm. Therefore, as $1<m\gamma<p^*$, there is $u\in E$ and a subsequence
(still denoted by $u_n$) such that $u_n$ converges weakly to $u$ in $E$, $\|u_n\|_{L^{m\gamma}(\O)}\to \|u\|_{L^{m\gamma}(\O)}$, and
\begin{equation*}
\lambda_M=\liminf_{n\to +\infty} \|u\|^{m\gamma}\geq \|u\|^{m\gamma}.
\end{equation*}
 Consequently, $\|u\|_{L^{m\gamma}(\O)}=1$ and so $\|u\|^{m\gamma}\geq\lambda_M$ which implies that $\|u\|^{m\gamma}=\lambda_M$.
\\
Moreover, there exists a Lagrange multiplier $\mu$ such that

$$mM(\|u\|^m)\int_\O|\textit{D}_r u|^{m-2}\textit{D}_r u\textit{D}_r v  =\mu m\gamma\int_\Omega|u|^{m\gamma-1} v,\; \forall v\in E.$$
If $v=u$ we have $m\lambda_M=\mu m\gamma$, that is $\lambda_M=\mu m\gamma$ and so
\begin{eqnarray*}
\begin{cases}
M(\|u\|^m)\Delta^r_m u =\lambda_M|u|^{m\gamma-1}  &\mbox{in}\quad \Omega, \\
u=\left(\frac{\partial}{\partial \nu}\right)^k u=0, \quad &\mbox{on}\quad \partial\Omega, \quad k=1, 2,.....  , r-1.
\end{cases}
\end{eqnarray*}
$(b)$ In general $\lambda_M$ is not attained if $M$ satisfies $(M_3)$. The typical example is $\widehat{M}(t)=t^\gamma H(t)$, where
$H:\R_+\to \R_+$ is strictly decreasing
with $\lim_{t\to \infty}H(t)=l>0$
  and $t^\gamma H(t)$ is strictly increasing on $\R_+$. So, $M$ is positive and satisfies $(M_3)$.
\\
Set
 $$E(u)=\dfrac{  \widehat{M}(\|u\|^m) }{\int_\O|u|^{m\gamma} dx}.$$
The monotonicity of $H$ involving $E(\alpha u)<E(u)$ \;for all\; $\alpha>1,$ and $u\neq0.$ It means
clearly that $\lambda_M$ is not attained. More exactly, one can
find easily examples of $H$ such as $H(t)=1+\beta(t+1)^{-1}$ (with $\beta>0$) and which also satisfy $(M_1)$.\qed
\subsection{ \bf Proof of  Theorem \ref{th24}}
First of all observe that $(M_1)$ (with $\tau_0=0$) implies
that for each $\tau_1>0$, we have
\begin{eqnarray}\label{MMM}
 \frac{\widehat{M}(\tau)} {\tau^\gamma} \leq \frac{ \widehat{M}(\tau_1)}{\tau_1^\gamma} , \; \forall \tau \geq \tau_1.
\end{eqnarray}
  To prove  Theorem \ref{th24}, we shall verify the validity of the conditions of the standard mountain pass theorem \cite{R}. Since $(M_3)$ implies $(M_2)$, Proposition \ref{th13} holds. Consequently,  $I$ satisfies the $(PS)$ condition if $m\geq 2$ (respectively the $(C)$ condition if $1<m<2$). By combining $(H_{2})$ and $(H'_{3})$ (at $0$), we can find $\epsilon_{0} > 0$ small enough and $C_0>0$ such that
$F(x ,s)\leq(\frac{\lambda_{M} }{m}-\epsilon_{0}) | s |^{m\gamma}+ C_0 | s |^{p^*} \mbox{ for all } (x,s)\in \Omega \times\mathbb{R}.$
 Also recall that $(M_3)$ implies $(i)$ of Lemma \ref{eig.posi} which with \eqref{Gagliardo} implies
 \begin{align*}
 I(u)  &\geq  \frac{1}{m}\widehat{M}(\|u\|^m)  -  (\frac{\lambda_{M} }{m}-\epsilon_{0}) \int_\Omega |u|^{m\gamma} - C_0 \int_\Omega |u|^{p^*}\\&\geq \frac{\epsilon_{0}}{\lambda_{M}}\widehat{M}(\|u\|^m)  - C'_0 || u ||^{p^*},\; C'_0>0.
 \end{align*}
Set $ \|u\|=\rho$ with $0<\rho \leq 1$, thus using $(M_3)$, we deduce
\begin{align*}
 I(u)  \geq \frac{C\epsilon_{0}}{\lambda_{M}}\rho^{m\gamma} - C'_0 \rho^{p^*}\geq \rho^{m\gamma}(\frac{C\epsilon_{0}}{\lambda_{M}} - C'_0 \rho^{p^*-m\gamma}).
 \end{align*}
     Choose $\rho =\inf(1,(\frac{C\epsilon_{0}}{2C'_0\lambda_{M}})^{\frac{1}{p^*-m\gamma}})$ and $\alpha=\frac{C\epsilon_{0}}{2\lambda_{M}}\rho^{m\gamma} >0$, then we have $I(u)\geq \alpha$ for all $\|u\|=\rho$.

 On the other hand, using $(H'_{3})$ (at infinity) and part (i) of Lemma \ref{lam pos}, then for $\epsilon_0>0$ small enough, we can find a positive constant $C_0$ and $\varphi \in E\setminus\{0\}$ such that
\begin{eqnarray}\label{e44}
|F(x,s)| \geq( \frac{\lambda_{M}}{m} +2\epsilon_{0}) | s |^{m\gamma}-C_0,\;\forall (x,s) \in \O \times\mathbb{R},
\end{eqnarray}
and
\begin{eqnarray}\label{e44'}
\lambda_{M}\int_\O|\varphi|^{m\gamma} \leq  \widehat{M}(\|\varphi\|^m)\leq (\lambda_{M}+m\epsilon_{0})\int_\O|\varphi|^{m\gamma}.
\end{eqnarray}
 Set $v=t\varphi,\;t\geq 1$ and using \eqref{e44}, we obtain
\begin{eqnarray}\label{e33}
I(v)&\leq& \frac{1}{m}\widehat{M}(t^m\|\varphi\|^m)
 -( \frac{\lambda_{M}}{m} +2\epsilon_{0})t^{m\gamma} \int_\Omega |\varphi |^{m\gamma}+C_0 |\Omega|\nonumber\\
&\leq& \frac{1}{m}\left(\dfrac{\widehat{M}(t^m\|\varphi\|^m)}{t^{m\gamma}}
-(\lambda_{M}+m\epsilon_{0})
\displaystyle\int_\O|\varphi|^{m\gamma}\right)t^{m\gamma} -\epsilon_0t^{m\gamma}\displaystyle\int_\O|\varphi|^{m\gamma}+C_0 |\Omega|.
\end{eqnarray}
     Using now \eqref{MMM} with $\tau_1=\|\varphi\|^m$, we obtain \begin{eqnarray}\label{ttm}
  \dfrac{\widehat{M}(t^m\|\varphi\|^m)}{ t^{m\gamma}}\leq \widehat{M}(\|\varphi\|^m),\; \mbox{ for all } t\geq 1.\end{eqnarray}
   So, from \eqref{e44'} and \eqref{e33}, we derive
  \begin{eqnarray*}\label{e333}
I(t\varphi) \leq \frac{1}{m}\left( \widehat{M}(\|\varphi\|^m) -(\lambda_{M}+m\epsilon_{0})\displaystyle\int_\O|\varphi|^{m\gamma}\right)t^{m\gamma} -\epsilon_0t^{m\gamma}\displaystyle\int_\O|\varphi|^{m\gamma}+C_0 |\Omega| \leq   -\epsilon_0t^{m\gamma}\displaystyle\int_\O|\varphi|^{m\gamma}+C_0 |\Omega|.\nonumber
\end{eqnarray*}
Choose $t$ large enough, we deduce that $I(v)<0$. In conclusion, $I$ satisfies the mountain pass geometry which ends the proof of Theorem \ref{th24}. \qed
\subsection{\bf Proof of Proposition \ref{sublinear}}
In view of $(H'_{1})$ and for $\epsilon_{0} > 0$ small enough, there exists $C_0>0$ such that
$F(x ,s)\leq(\frac{\lambda_{M} }{m}-\epsilon_{0}) | s |^{m\gamma}+ C_0$  for all $(x,s)\in \Omega \times\mathbb{R}.$ According to $(M_3)$ we derive
$$ I(u)  \geq \frac{1}{m}\widehat{M}(\| u\|^{m})-(\frac{\lambda_{M} }{m}-\epsilon_{0}) \int_\Omega |u|^{m\gamma} -C_0|\O|\geq\frac{\epsilon_{0}}{\lambda_{M}}\widehat{M}(\| u\|^{m})  - C_0|\O| \geq C\| u\|^{m\gamma}- C_0|\O|.
$$
Consequently, $I(u)\to \infty$ as $\|u\|\to \infty$. Then any $(PS)$ sequence is bounded. Since $(H'_1)$ implies $(H_2)$ as $ 1<m\gamma <p^*$, then from the proof of step $2$ of Proposition \ref{th13} we deduce that any bounded $(PS)$ sequence  has a strong convergent subsequence. So $I$ verifies the $(PS)$ condition.\qed
\\
\subsection{\bf Proof of Theorem \ref{thm1.6}}
Recall first some basic notations of Krasnoselskii's genus, which can be found in \cite{R}. Let $E$ be a Banach space and $$\Sigma=\{A \subset X -\{0\}: A \mbox{ is closed in } X \mbox{ and symmetric with respect to }0\}.$$
\begin{definition}(See \cite{R}) For $A \in \Sigma$, we say genus of $A$ is $n$ denoted by $\gamma(A) = n$
if there is an odd map $\phi \in C(A, \R^n \backslash\{0\})$ and $n$ is the smallest integer with
this property.
\end{definition}
We invoke the following abstract theorem  based on the Krasnoselskii genus theory to prove Theorem \ref{thm1.6}.
\begin{theorem}\label{Lemma2.3}(See \cite{R})
 Let $I \in C^1(E)$ be an even functional satisfying the $(PS)$
condition. For $c\in \R$ and each $n \in \N^*$, set $K_c=\{u\in E:\; I'(u)=0, \;I(u)=c\}$ and
$\Sigma_n=\{A\in \Sigma:\gamma(A)\geq n\},\;\;c_n=\inf_{A\in \Sigma_n}\sup_{u\in A}I(u).$
\begin{enumerate}
  \item [$(i)$] If  $\Sigma_n\neq \emptyset$ and $-\infty<c_n<0$, then $c_n$ is a critical value of $I$.
  \item [$(ii)$] If there exists $\varrho\in \N$  such that $c_n=c_{n+1}=...=c_{n+\varrho}=c\in \R,$ and $c\neq I(0)$, then $\gamma(K_c)\geq \varrho+1.$
\end{enumerate}
\end{theorem}
We will verify the conditions of Theorem \ref{Lemma2.3}.\\
 Clearly, $I$ is even with $I(0) = 0$ and from Proposition \ref{sublinear} $I $ satisfies the $(PS)$ condition.  We shall now prove that for any $n\geq 2, \; n\in \N$, one has $\Sigma_n \neq  \emptyset$. In fact, we can find $\phi_1, \phi_2,..., \phi_n\in C_c^{2r}(\O)$  satisfying $\|\phi_i\|_{L^2(\O)}=1$ and $supp(\phi_i)\cap supp(\phi_j)=\emptyset$ if $i\neq j,\; 1\leq i,\;j\leq n$. Set
$$E_n= span\{ \phi_1, \phi _2,.., \phi_n\}=\{u=\sum_{i=1}^{n}\lambda_i \phi_i,\; \lambda_i \in R\} \subset E \cap L^2(\O),$$  and for $0<\sigma<1$  $$ S^\sigma_n=\{u\in E_n : \|u\|_{L^2(\O)}=\sigma\}=\{u\in E_n,\;
\sum_{i=1}^{n}\lambda_i^2=\sigma^2\}.$$
Consider the map $h:S_n^\sigma \to S^{n-1}$  defined by:
$$h(u)=(\frac{\lambda_1}{ \sigma},\;\frac{\lambda_2}{\sigma},...,\;\frac{\lambda_n}{\sigma}), \; u \in S_n^\sigma,$$
 where $S^{n-1}=\{(\beta_1, \beta_2,...,\beta_n)\in \R^n:\sum_{i=1}^{n}\beta_i^2=1\}$ is the sphere of dimension $n-1$. Clearly $h$ is an  homeomorphic odd map, which means that $\gamma(S_n^\sigma)=n$ (see \cite{R}), then $S_n^\sigma \in \Sigma_n$ and the claim is well proved.\\
Therefore, $c_n$ is well defined and since point $1$ of Proposition \ref{sublinear} implies that $I$ is bounded from below, so $-\infty <c_n$.

In order to apply Theorem \ref{Lemma2.3} we have to prove that
$c_n<0.$ Fix $n\in \N^*$. Indeed, since all norms in $E_n$ are equivalent, there exists
 $C_n > 0$ such that
\begin{eqnarray}\label{equiv}
\frac1{C_n} \|u\|\leq \|u\|_{L^2(\O)}\leq C_n \|u\|_{L^{m\gamma}(\O)}.
 \end{eqnarray}
According to $(H'_2)$ one has for every $A>0$ there is $s_A>0$ such that
\begin{eqnarray}\label{negF}
F(x,s)\geq A| s|^{m\gamma}, \mbox{ for all } x\in \O,\; |s|\leq s_A.
\end{eqnarray}
Set $M_n=\max \{\|\phi_i\|_{L^\infty(\O)},\; 1\leq i\leq n\}$, then for  $\sigma =\inf(\frac 12,\frac{ s_A}{2nM_n})$ we have $\|u\|_{L^\infty(\O)}\leq \frac{ s_A}{2},\; \forall u\in S^\sigma_n .$
Choose now $A=\beta C^{2m\gamma}_n$, so by combining  $(M'_2)$, \eqref{equiv} and \eqref{negF} we deduce that
\begin{eqnarray*}
I( u)=\frac{1}{m}\widehat{M}(\| u\|^{m}) -\int_\O F(x, u)
\leq\frac{\beta}{mC_n^{m\gamma}}\left(C^{2m\gamma}_n\| u\|_{L^2(\O)}^{m\gamma} -\frac{mA}{\beta}\|u\|_{L^2(\O)}^{m\gamma}\right)
\leq\frac{\beta(1-m)C_n^{m\gamma}}{m}\sigma^{m\gamma}<0,\;\mbox{ for all } u\in S^\sigma_n.
\end{eqnarray*}
As $m>1$, then $\frac{\beta(1-m)C_n^{m\gamma}}{m}\sigma^{m\gamma}<0$, thus $\sup_{u\in S^\sigma_n }I(u)<0$ and $c_n=\inf_{A\in \Sigma_n}\sup_{u\in A}I(u)<0$.
In conclusion, from point $1$ of Theorem \ref{Lemma2.3} we derive that $c_n$ is a critical value of $I$ and since $n$ is arbitrary, then point $2$ of Proposition \ref{sublinear} implies that $I$ admits infinitely many nontrivial critical points. The proof
is completed.\qed

\subsection{\bf Proof of Theorem \ref{th14}}
In view of Proposition \ref{sublinear}, $I$ satisfies the $(PS)$ condition and $I$ is bounded from below. Also, we claim that  $c=\inf\{I(u),\; u\in E\}<0$. Indeed, as $f$ verifies $(H''_2)$ at $0$, so for $\epsilon_{0} > 0$ small enough there exists $s_0> 0$ such that
\begin{eqnarray}\label{ee44}
 F(x ,s)\geq (\frac{\lambda_{M} }{m}+2\epsilon_{0}) | s |^{m\gamma}  \mbox{ for all } (x,s)\in \Omega \times [-s_0,\;s_0].
 \end{eqnarray}
Taking into account that $\overline{C_c^r(\O)}=E$, and according to Lemma \ref{lam pos}, there is
  $\phi \in C_c^r(\O)\setminus \{0\}$ such that
 \begin{eqnarray}\label{ee44'} \widehat{M}(\|\phi\|^{m})\leq (\lambda_{M}+m\epsilon_{0})\int_\O|\phi|^{m\gamma}.
\end{eqnarray}
  Using now \eqref{MMM} (which is a consequence of $(M_3)$), we derive
\begin{eqnarray}\label{nn}
  \dfrac{\widehat{M}(t^m\|\phi\|^m)}{ t^{m\gamma}}\leq \widehat{M}(\|\phi\|^m),\; \forall t\leq 1.
  \end{eqnarray}
 Set $t=\inf(\frac{s_0}{\|\phi\|_{L^{\infty}(\O)}}, \frac 12)$. So $t\phi \in [-s_0,s_0]$ and by combining \eqref{ee44}-\eqref{nn}, we get
 \begin{eqnarray}
I(t\phi)\leq \frac{t^{m\gamma}}{m}\left(\widehat{M}(\|\phi\|^{m})-(\lambda_{M}+m\epsilon_{0})\int_\O|\phi|^{m\gamma}\right)- \epsilon_{0}t^{m\gamma}\int_\O|\phi|^{m\gamma} <0.
\end{eqnarray}

Therefore, $c=\inf\{ I(u),\; u \in E\}<0$, invoking Ekeland's variational principle, we deduce that $c$ is a nontrivial critical value which achieves the proof of Theorem \ref{th14}.\qed
\section{Proof of Theorem \ref{th46}}\label{section 8}
For the convenience of the reader, we recall in this section some theories on spaces $L^{m(x)}(\O)$ and $W^{r,m(x)}(\O)$ which
we call generalized Lebesgue-Sobolev spaces. Denote
$$C_+(\overline{\O})=\left\{m(x);\; m(x)\in C(\overline{\O}),\; m(x)>1, \;\forall \; x\in \overline{\O}\right\}.$$
For any $m\in C_+(\overline{\O})$, we introduce the variable exponent Lebesgue space
$$L^{m(.)}(\O)=\left\{u: u \mbox{ is a measurable real-valued function such that } \int_\O |u(x)|^{m(x)}dx<\infty\right\},$$
endowed with the so-called Luxemburg norm
$$\|u\|_{L^{m(x)}(\O)}=|u|_{m(.)}=\inf \left\{\mu>0;\int_\O\left|\frac{u(x)}{\mu}\right|^{m(x)}dx\leq 1\right\},$$
which is a separable and reflexive Banach space. A thorough variational analysis of the problems with variable exponents has been developed in
the recent monograph by R$\breve{a}$dulescu and Repov$\breve{s}$ \cite{RD} (we refer also the reader to \cite{Diening,FF,KR,Y}).

\begin{prop}[see \cite{Y}]
 The space $(L^{p(x)}(\O), |.|_{p(x)})$ is separable, uniformly convex, reflexive and its conjugate space is $L^{p(x)}(\O), |.|_{q(x)}$ where $q(x)$ is the conjugate function of $p(x)$ i.e
$$\frac{1}{p(x)}+\frac{1}{q(x)}=1,\;\;\forall x\in \O.$$
For all $u\in L^{p(x)}(\O)$ and $v\in L^{q(x)}(\O)$ the H\"older's type inequality
\begin{eqnarray*}
  \left|\int_\O uvdx\right|\leq\left(\frac{1}{p^-}+\frac{1}{q^-}\right)|u|_{p(x)}|v|_{q(x)} \end{eqnarray*}
  holds true.
\end{prop}

The inclusion between Lebesgue spaces also generalizes the classical framework,
namely if $0<|\O|<\infty$ and $p_1$, $p_2$ are variable exponents such that $p_1 \leq p_2$ in $\O$ then there exists a continuous embedding $L^{p_2(x)}(\O)\to L^{p_1(x)}(\O)$. An important role in manipulating the generalized Lebesgue-Sobolev spaces is played by the $m(.)-$modular of the $L^{m(.)}(\O)$ space, which is the modular $\rho_{m(\cdot)}$ of the space $L^{m(\cdot)}(\O)$
\begin{equation*}
   \rho_{m(\cdot)}(u):=   \int_{\Omega} |u|^{m(x)} \,dx.
\end{equation*}

\begin{lem}\label{lemmaineq}
If $u_n, u \in L^{m(\cdot)}$ and $m_{+} < +\infty$, then the following
properties hold:
\begin{enumerate}
\item $|u|_{m(\cdot)} > 1 \Rightarrow
|u|_{m(\cdot)}^{m_{-}} \leq \rho_{m(\cdot)}(u) \leq |u|_{m(\cdot)}^{m_{+}}$;

\item $|u|_{m(\cdot)} < 1 \Rightarrow   |u|_{m(\cdot)}^{m_{+}}
\leq \rho_{m(\cdot)}(u) \leq |u|_{m(\cdot)}^{m_{-}}$;

\item $|u|_{m(\cdot)} < 1$ (respectively $= 1; > 1) \Longleftrightarrow
 \rho_{m(\cdot)}(u) < 1$ (respectively $= 1; > 1$);\label{gggg}

\item $|u_n|_{m(\cdot)} \to 0$ (respectively
$\to +\infty) \Longleftrightarrow \rho_{m(\cdot)}(u_n) \to 0$
 (respectively $\to +\infty$);
\item $\lim_{n\to \infty}|u_n-u|_{m(x)}=0 \Longleftrightarrow \lim_{n\to \infty}\rho_{m(\cdot)}(u_n-u)=0$.
\end{enumerate}
\end{lem}

The Sobolev space with variable exponent $W^{r,m(x)}(\O)$ is defined as
\begin{equation*}
   W^{r,m(x)}(\Omega):=\Big\{u \in L^{m(x)}(\Omega):
D^\alpha u \in L^{m(x)}(\Omega),\;|\alpha|\leq k \Big\},
\end{equation*}
where $D^\alpha u=\frac{\partial^{|\alpha|}}{\partial x_1^{\alpha_1}\partial x_2^{\alpha_2}...\partial x_N^{\alpha_N}u}$,
 with $\alpha=(\alpha_1,...\alpha_N)$ is a multi-
 index and $|\alpha|=\sum_{i=1}^{N}\alpha_i.$ The space $ W^{r,m(x)}(\Omega)$ is a reflexive and separable Banach space if $1<m_-\leq m_+<+\infty$ and equipped with the norm
$$\|u\|_{r,m(x)}:=\sum_{|\alpha\leq k}|D^\alpha u|_{m(x)}.$$

Let $W_0^{r,m(x)}(\O)$ denote the completion of $C_0^\infty (\O)$ in $W^{r,m(x)}(\O)$. As shown in (\cite{Diening},
Corollary $11.2.4$), the space $W_0^{r,m(x)}(\O)$ coincides with the closure in $W^{r,m(x)}(\O)$ of the set of all $W^{r,m(x)}(\O)$-functions with compact support.
\begin{prop}[see \cite{FZ1}]\label{embedding} For $m,p\in C_+(\overline{\O})$ such that $p(x)\leq p^*(x)$ for all $x\in \overline{\O}$, there is a continuous embedding $W^{r,m(x)}(\O)\hookrightarrow L^{p(x)}(\O)$. If we replace $\leq$ with $<$, the embedding is compact.
\end{prop}
In order to prove Theorem \ref{th46}, we need to use again the symmetric mountain pass theorem. We first have the following lemmas.
\begin{lem}\label{lem1}
Assume that $f$ and $M$ verify $(H^{exp}_1)$-$(H^{exp}_2)$ and $(M_1)$-$(M_2)$ respectively, then
\begin{itemize}
  \item[$1)$] $I$ satisfies $(PS)$ condition for\; $m_+\geq m_-\geq2$.
  \item[$2)$] $I$ satisfies $(C)$ condition for\; $2>m_+\geq m_->1$.
\end{itemize}
\end{lem}
{\bf Proof.} According to $(M_1)$ and $(H^{exp}_1)$, there exists $C_0>0$ such that
\begin{eqnarray}\label{er1}
C|f(x,s)|^{q^*(x)}- C_0\leq s f(x,s)- m_+\gamma F(x,s),\; \forall (x,s)\in \Omega\times \R
\end{eqnarray}
and
\begin{eqnarray}\label{er2}
-C_0 \leq \gamma \widehat{M}(\tau)-\tau M(\tau), \; \forall \tau\geq 0.
\end{eqnarray}
Since we assume that $N> rm(x)$, we may easily see that
 \begin{eqnarray}\label{err}
 m_+\geq m_-> \frac{1}{(q^*)_-}+1 \mbox{ if  } m_+\geq m_-\geq2.
\end{eqnarray}
Let $\{u_n\}$ be a $(PS)$ sequence of $I$ if $m_+\geq m_-\geq2$ ({ respectively} $(C)$ sequence if $1<m_-\leq m_+<2$). Two possible cases arise: either $\{u_n\}$ admits a subsequence which converges strongly to $0$ in $W_0^{r,m(x)}(\O)$ and so we have done, or there exist $\eta_0>0$ and $n_0 \in \N$ such that $\|u_n\| \geq \eta_0$ for all $n\geq n_0$. Observe that, by Lemma \ref{lemmaineq} we get for all $n \in \N$,
$$
\int_\O \frac{|D_r u_n|^{m(x)}}{m(x)}dx \geq \frac{1}{m_+}\min\left\{\|u_n\|^{m_+},\|u_n\|^{m_-}\right\} \geq \frac{1}{m_+}\min\left\{\eta_0^{m_+}, \eta_0^{m_-}\right\}.
$$
So according to $(M_2)$, there is ${m}_{\eta_0}>0$ such that
 \begin{eqnarray}\label{ing1}
 M\left(\int_\O \frac{|D_r u_n|^{m(x)}}{m(x)}dx\right)\geq {m}_{\eta_0}>0, \forall n\geq n_0.
\end{eqnarray}
\\
\textbf{Step 1.} $\{u_n\}$ is bounded in $W_0^{r,m(x)}(\O)$. Taking into account \eqref{kir deriv} and  \eqref{ing1}, we can deduce that
\begin{align*}
m_{\eta_0}\rho_{m(x)}(\textit{D}_ru_n) \leq M\left( \int_\O \frac{|D_r u_n|^{m(x)}}{m(x)}dx\right)\int_\O|D_r u_n|^{m(x)}dx =&\langle I'(u_n),u_n\rangle+\int_{\O} f(x,u_n)u_ndx.
\end{align*}

Applying H\"older's inequality to the second term in the right-hand side and using Proposition \ref{embedding}, we obtain
\begin{align}\label{kirch hk}
m_{\eta_0}\min\left\{\|u_n\|^{m_-}, \|u_n\|^{m_+}\right\} \leq m_{\eta_0} \rho_{m(x)}(\textit{D}_ru_n) \leq \langle I'(u_n),u_n\rangle&+ C\left|f(x,u_n)\right|_{q^*(x)}\|u_n\|.
\end{align}
From \eqref{kir func}, \eqref{kir deriv} and \eqref{er1}, \eqref{er2}, one has
\begin{align}\label{ing2}
m_+\gamma I(u_n)-\langle I'(u_n), {u_n}\rangle&=\left[m_+\gamma
\widehat{M}\left( \int_\O \frac{|D_r u_n|^{m(x)}}{m(x)}dx\right)-M\left( \int_\O \frac{|D_r u_n|^{m(x)}}{m(x)}dx\right)\int_\O |D_r u_n|^{m(x)}dx\right]\nonumber\\
&\quad + \int_\O\left[f(x,u_n)u_n-m_+\gamma F(x,u_n)\right]dx\nonumber\\
& \geq m_+\left[\gamma\widehat{M} \left( \int_\O \frac{|D_r u_n|^{m(x)}}{m(x)}dx\right)-M\left( \int_\O \frac{|D_r u_n|^{m(x)}}{m(x)}dx\right)\int_\O \frac{|D_r u_n|^{m(x)}}{m(x)}dx\right]\nonumber\\
& \quad + C\int_\O|f(x,u_n)|^{q^*(x)}dx- C_0\nonumber\\
&\geq -C_0(m_++1)+ C\int_\O|f(x,u_n)|^{q^*(x)}dx.
 \end{align}

Taking into account that $\{u_n\}$ is a $(PS)$ sequence if $m_+\geq m_-\geq2$ (respectively $(C)$ sequence if $1<m_-\leq m_+<2$), then we deduce from \eqref{ing2} that
\begin{equation}\label{lmm"}\int_{\Omega} | f(x , u_{n}) |^{q^*(x)}dx\leq C(1+\|u_{n}\|), \left(\mbox{ respectively }\int_{\Omega} | f(x , u_{n})|^{q^*(x)}dx\leq  C\right).
  \end{equation}
Combining now \eqref{kirch hk} with \eqref{lmm"}, it follows that
\begin{equation*}
m_{\eta_0}\min\left\{\|u_n\|^{m_-}, \|u_n\|^{m_+}\right\}\leq
\begin{cases}
  C\left(1 + || u_{n} ||^{{\frac{1}{(q^*)_-}+1}}\right), & \mbox{if } m_+\geq m_-\geq2, \\
  C(1 + \| u_{n} \| ), & \mbox{if } 1<m_-\leq m_+<2.
\end{cases}
\end{equation*}

  Clearly, the $(C)$ sequence $u_n$ is bounded in $W_0^{r,m(x)}(\O)$ if $1<m_-\leq m_+<2$ and thanks to \eqref{err} the $(PS)$ sequence is also bounded if $m_+\geq m_-\geq2$.
 \\\\
{\bf Step 2.} We shall prove that the bounded sequence $\{u_n\}$ has
a strong convergent subsequence in $E$. Indeed, there exist a subsequence (denoted by $u_n$) and $u \in W_0^{r,m(x)}(\O)$ such that $u_n$ converges to
$u$ weakly in $W_0^{r,m(x)}(\O)$ and strongly in $L^1(\O)$. Also $u_n$ and $u_n - u$ are bounded in $L^{p^*(x)}(\O)$ and from \eqref{kir deriv}, we get
\begin{eqnarray}\label{e:1221zz}
M\left( \int_\O \frac{|D_r u|^{m(x)}}{m(x)}dx\right)\int_\O|\textit{D}_r u_n|^{m(x)-2}\textit{D}_r u_n\textit{D}_r (u_n-u)dx=\langle I'(u_n), u_n-u\rangle+\int_\Omega f(x,u_n)(u_n-u)dx.\nonumber\\
\end{eqnarray}
 By the virtue of the  condition $(H^{exp}_2)$, one has for every $\epsilon \in (0,1)$,  there exists $C_{\epsilon}>0$ such that
 \begin{equation}
\label{e:1221}
|f(x,s)|\leq \epsilon |s|^{p^*(x)-1}+C_{\epsilon},\,\, \forall (x,s)\in \overline\O\times \R.
\end{equation}

As $u_n$ converges strongly in $L^1(\O)$, there exists $N_\epsilon$ such that $ \int_{\O}|u_n-u|\leq \frac{\epsilon}{C_\epsilon},\; \forall n\geq N_\epsilon$. So, in view of \eqref{e:1221} and H\"older's inequality, we obtain
\begin{align*}\label{K3.6}
\left|\int_\Omega f(x,u_n)(u_n-u)dx\right|&\leq \epsilon\int_{\O}|u_n|^{p^*(x)-1}|u_n-u|dx + C_\epsilon \int_{\O}|u_n-u|dx\nonumber \\
&\leq \epsilon\left||u_n|^{p^*(x)-1}\right|_{\frac{p^*(x)}{p^*(x)-1}}|u_n-u|_{p^*(x)} + C_\epsilon \int_{\O}|u_n-u|dx\nonumber\\
&\leq \epsilon \max\left\{|u_n|^{(p^*)_+-1}_{p^*(x)},|u_n|^{(p^*)_--1}_{p^*(x)}\right\}|u_n-u|_{p^*(x)} + C_\epsilon \int_{\O}|u_n-u|dx\nonumber\\
&\leq C\epsilon, \;\forall n\geq N_\epsilon.
\end{align*}

Consequently, $\int_\Omega f(x,u_n)(u_n-u)dx$ converges to $0$, and since $I'(u_n) \to 0$ in $[W_0^{r,{m(x)}}(\O)]'$ and $(u_n-u)$ is bounded in $W_0^{r,m(x)}(\O)$, we deduce from \eqref{e:1221zz} that
$$M\left( \int_\O \frac{|D_r u|^{m(x)}}{m(x)}dx\right)\int_\O|\textit{D}_r u_n|^{m(x)-2}\textit{D}_r u_n\textit{D}_r (u_n-u)dx\to 0.$$

Using again \eqref{ing1}, yields  $\int_\O|\textit{D}_r u_n|^{m(x)-2}\textit{D}_r u_n\textit{D}_r (u_n-u)dx \mbox{ converges to } 0.$
Arguing now as in \cite{CP}, we obtain $\limsup_{n\to +\infty}\rho_{m(x)}(\textit{D}_ru_n-\textit{D}_r u)=0$, which implies that $\{u_n\}$ converges strongly to $u$ in $W_0^{r,{m(x)}}(\O)$.\qed
\begin{lem}\label{lem2}
Suppose that $(H^{exp}_3)$ and $(M_1)$ are satisfied. Then for any finite dimensional subspace $W\subset W_0^{r,{m(x)}}(\O)$ there is $R= R(W)>0$ such that $I(u)\leq 0$ for $u\in W$, $\|u\|\geq R$.
\end{lem}
{\bf Proof.} In view of $(H^{exp}_3)$ we have for all $A>0$, there is $C_A>0$ such that
$$F(x,s)\geq A|s|^{m_+\gamma}- C_A,\; \mbox{ for all } (x,s)\in \O\times \mathbb{R}.$$
Using again $(M_1)$ we derive that
$$\widehat{M}(\tau)\leq C_1\tau^{\gamma}-C_2,\; \forall \tau \geq 0,\mbox{  where } C_1 =\dfrac{\widehat{M}(\tau_0)}{\tau_0^{\gamma}} \mbox{ and }C_2>0.$$
Consequently, we obtain
\begin{eqnarray*}
 I(u) &=& \widehat{M}\left( \int_\O \frac{|D_r u|^{m(x)}}{m(x)}dx\right)-\int_\O F(x,u)dx\\
 &\leq& \frac{C_1}{m_-^{\gamma}}[\rho_{m(x)}(\textit{D}_ru)]^\gamma-A\|u\|^{m_+\gamma}_{L^{m_+\gamma}(\O)} +C_A-C_2 \\
 &\leq& \frac{C_1}{m_-^{\gamma}}\max\left\{\|u\|^{m_-\gamma },\|u\|^{m_+\gamma}\right\}-A\|u\|^{m_+\gamma}_{L^{m_+\gamma}(\O)} +C_A-C_2.
\end{eqnarray*}

 Let $W$ be a fixed finite dimensional subspace of $W_0^{r,{m(x)}}(\O)$, as $\|.\|$ and
$\|u\|_{L^{m_+\gamma}(\O)} $ are equivalent norms on $W$. Take $R= R(W)>0$, then  for all $u\in W$ with $\|u\|\geq R$ we obtain
 $$ I(u) \leq \|u\|^{m_+\gamma}\left(\frac{C_1}{m_-^{\gamma}}-AC_W\right) +C_A-C_2,$$
where $C_W>0$ is such that $\|u\|_{L^{m_+\gamma}(\O)}\geq \|u\|$. Choosing $A=\dfrac{2C_1}{C_Wm_-^{\gamma}}$ in the last inequality, so we derive  $I(u)< 0$ for all $u\in W$ and $\|u\|\geq R$. This completes the
 proof.\qed
 \\

According to \cite{FHHSPZ,Le,MMO,YZS}, we notice that $W_0^{r,m^+}(\O)\subset W_0^{r,m(x)}(\O)$. Consider $\{e_1, e_2, . . .\}$, a Schauder basis of the space $W_0^{r,m^+}(\O)$ which means that each $x\in W_0^{r,m^+}(\O)$ has a
unique representation $x=\sum_{i=1}^{\infty}a_ie_i$, where $a_i$ are real numbers. We consider for each $j\in \N^*$ that $E_j$, the subspace of $W_0^{r,m^+}(\O)$ generated by $j$ vectors $\{e_1, e_2,..,e_j\}$.
Clearly $E_j$ is subspace of $W_0^{r,m(x)}(\O)$.
 \\

As a consequence, the linear projection onto $E_j$ i.e.,
$P_j: W_0^{r,m^+}(\O)\to E_j,\; P_j(x)= \sum_{i=1}^{j}a_ie_i$ is continuous for all $j\in \N^*$ (see \cite{FHHSPZ}) \footnote{$P_j$ are uniformly bounded, that is there exists $C>0$ such that $\|P_j(x)\|\leq C\|x\| $ for each $j\in \N^*$ and all $x\in W_0^{r,m^+}(\O)$}.
  Therefore, $F_j= N( P_j)$ (the kernel of $P_j$) is a topological complement of $E_j$, that is $E_j\oplus F_j=W_0^{r,m^+}(\O)$. Fix $\rho\geq 0$ and set $S_j^\perp(\rho)=\left\{u\in F_j \mbox{ such that }\|u\|=\rho\right\}$. Then we have the following lemma.

\begin{lem}\label{lem3}
Suppose that $(H^{exp}_2)$ is satisfied. Then for all $\rho\geq 0$, there exist $j_0$ and $\alpha>0$ such that $I(u)\geq \alpha$, $\forall u\in S_{j_0}^\bot(\rho)$.
\end{lem}
{\bf Proof.} The idea consists in applying the Schauder basis of $W_0^{r,m(x)}(\O)$. For this purpose, we proceed by the following Claims.
\\
{\bf Step $1$.} Suppose that $(H^{exp}_2)$ is satisfied and let
$$\beta_j:=\sup_{S_j^\perp(\rho)} \int_{\Omega} |F(x,u)|dx.$$
Then $\beta_j \rightarrow 0, \mbox{ as } j\rightarrow \infty.$
\\

We argue by contradiction. Suppose that there exist $m_0>0$ and a subsequence (denoted by $\beta_j$) such that
$m_0<\beta_j,\forall j\in \N^*.$
From the definition of $\beta_j$, there exists $u_j\in S_j^\perp(\rho)$ such that
\begin{eqnarray}\label{uuk}
m_0<\int_{\O}|F(x,u_j)|dx\leq \beta_j.
\end{eqnarray}

   As $\|u_j\|=\rho$,  $\{u_j\}$ is bounded in $L^{p^*(x)}(\O)$, consequently there exist a subsequence (denoted by $u_{j}$) and $u\in W_0^{r,m^+}(\O) $ such that $u_{j}$  converges weakly to $u$ and $\mbox{ a.e in }\; \Omega$. Fix $k\in \N^*$.\\ Since $F_k= N(P_k)$ is the kernel of $P_k$ and $P_k\circ P_{k+1}=P_k$, then $ F_{j}\subset F_k$ for all $j\geq k$. So as $u_{j} \in F_{j}$, then
  $P_k(u_{j})=0,$ for all $j\geq k$. Recall that $P_k$ is a linear continuous operator, then $P_k(u_{j})$
 converges weakly to $P_k(u)$ which implies that $P_k(u)=0$ and also $u=0$ as $u=\lim_{k\to
\infty}\sum_{i=1}^{k}a_ie_i=\lim_{k\to
\infty}P_k(u)$. Consequently, $u_j$ converges to $0$ $\mbox{ a.e in }\; \Omega$, and $F(x,u_{j})  \rightarrow 0 \mbox{ a.e in }\; \Omega$ as $F(x,0)=0$. By the virtue of $(H^{exp}_2)$, we have for every $\epsilon>0$ there is $C_\epsilon>0$ such that
 \begin{equation}\label{FF}
 |F(x,s) |\leq \epsilon | s |^{p^*(x)} + C_{\epsilon}, \;\;\mbox{ for all } (x,s)\in {\Omega}\times\mathbb{R}.
 \end{equation}Hence, for each measurable set $A\subset\O$ such that $|A|<\frac{\epsilon}{C_{\epsilon}}$ (where $|A|$ denotes the Lebesgue measure of $A$), we derive from \eqref{FF} the following
\begin{eqnarray*}
\int _{A} |F(x,u_{j})|dx \leq   \epsilon  \int _{A} | u_{j}|^{p^*(x)}dx
+ C_{\epsilon}|A| \leq C\epsilon.
\end{eqnarray*}
Taking into account that $\O$ is a bounded domain, then Vitali's theorem implies that $F(x,u_{j})\rightarrow 0$ in $L^{1}(\O)$ and in view of \eqref{uuk}, we obtain
$$
0< m_0\leq 0.
$$
 Thus, we reach a contradiction and and this concludes the proof of Step $1$.\qed
\\
{\bf Step $2$.} Let $\rho\geq0$, then for all $u\in S_j^\perp(\rho)$,
we have
\begin{eqnarray*}
I(u) &=& \widehat{M}\left( \int_\O \frac{|D_r u|^{m(x)}}{m(x)}dx\right)-\int_\O F(x,u)dx\\&\geq& \widehat{M}\left(\int_\O \frac{|D_r u|^{m(x)}}{m(x)}dx\right)-\beta_j.\end{eqnarray*}
As $\beta_j$ converges to $0$, we can choose $j=j_0$ large enough such that $$\beta_{j_0} \leq \frac {1}{2}\widehat{M}\left(\int_\O \frac{|D_r u|^{m(x)}}{m(x)}dx\right).$$
 Then,
 \begin{eqnarray*}
I(u)\geq \frac {1}{2}\widehat{M}\left(\int_\O \frac{|D_r u|^{m(x)}}{m(x)}dx\right)\geq \frac {m_0}{2m_+\gamma}\rho_{m(x)}(D_r u)\geq \frac {m_0}{2m_+\gamma}\min\{\|u\|^{m^-},\|u\|^{m^+}\}.
 \end{eqnarray*}
 Therefore, Lemma \ref{lem3} holds with $E^-=E_{j_0}$, $E^+=F_{j_0}$,  $\|u\|=\rho$, and $\alpha=\frac {m_0}{2m_+\gamma}\min\{\rho^{m^-},\rho^{m^+}\}$. Thus, Step $2$ is proved.\qed
\\{ \bf  Proof of Theorem \ref{th46}.}
By $(H^{exp}_4)$, we know that $I$ is even. Since $I(0)=0$, we conclude from Lemma \ref{lem1}, Lemma \ref{lem2}, Lemma \ref{lem3} and Theorem \ref{th45}, that problem \eqref{high} admits infinitely many distinct
pairs $(u_j,-u_j )_{j_\in \N^*}$, of critical points with unbounded energy.\qed
\subsection*{\bf Acknowledgments}
The first author was supported by the Tunisian Military Research Center for Science and Technology Laboratory LR19DN01. The second author would like to express his deepest gratitude to the International Centre for Theoretical Physics (ICTP), Trieste, Italie for providing him with an excellent atmosphere for doing this work.

Hamdani-Harrabi wish to thank Professors {\bf Dong Ye} and {\bf Nguyen Thanh Chung} for stimulating discussions on the subject.

\end{document}